\documentclass[11pt,onecolumn]{article}
\usepackage{multicol}
\usepackage[greek,english]{babel}
\usepackage[vcentering,dvips]{geometry}
\usepackage{amsfonts}
\usepackage{amssymb}
\usepackage{eucal}
\usepackage{amsxtra}
\usepackage{color}
\usepackage{setspace}
\usepackage{times}
\usepackage{titlesec}
\usepackage{sectsty}
\allsectionsfont{\large}
\usepackage{url}
\usepackage{doi}
\usepackage{hyperref}
\hypersetup{
    colorlinks,%
    citecolor=black,%
    filecolor=black,%
    linkcolor=black,%
    urlcolor=black
}

\addtocounter{section}{0} \numberwithin{equation}{section}
\newtheorem{theorem}{Theorem}[section]
\newtheorem{proposition}[theorem]{Proposition}
\newtheorem{lemma}[theorem]{Lemma}

\newcommand{\qed}{{$\hfill \Box$}}

\onehalfspace

\hypersetup{urlcolor=blue }

\begin{document}
\pagestyle{plain}

\title{\large
{\textbf{REAL HYPERSURFACES EQUIPPED WITH PSEUDO-PARALLEL STRUCTURE JACOBI OPERATOR IN $\mathbb{C}P^{2}$ AND
$\mathbb{C}H^{2}$}}}

\author{ \textbf{\normalsize{Konstantina Panagiotidou and Philippos J. Xenos}}\\
\small \emph{Mathematics Division-School of Technology, Aristotle University of Thessaloniki, Greece}\\
\small \emph{E-mail: kapanagi@gen.auth.gr, fxenos@gen.auth.gr}}
\date{}

\maketitle
\begin{flushleft}
\small {\textsc{Abstract}. Motivated by the work done in \cite{CK}, \cite{CHI}, \cite{NR2} and \cite{PS}, we classify real hypersurfaces in $\mathbb{C}P^{2}$ and $\mathbb{C}H^{2}$ equipped with pseudo-parallel structure Jacobi operator.}
\end{flushleft}
\begin{flushleft}
\small{\emph{Keywords}: Real hypersurface, Pseudo-parallel structure Jacobi operator, Complex projective space, Complex hyperbolic space.\\}
\end{flushleft}
\begin{flushleft}
\end{flushleft}

\section{Introduction}

A complex n-dimensional Kaehler manifold of constant holomorphic
 sectional curvature c is called a complex space form, which is
 denoted by $M_{n}(c)$. A complete and simply connected
 complex space form is complex analytically isometric to a complex
 projective space $\mathbb{C}P^{n}$, a complex Euclidean space
 $\mathbb{C}^{n}$ or a complex hyperbolic space $\mathbb{C}H^{n}$
 if $c>0, c=0$ or $c<0$ respectively.

Let $M$ be a real
 hypersurface in  a complex space form $M_{n}(c)$, $c\neq0$. Then an almost contact metric
 structure $(\varphi,\xi,\eta,g)$ can be defined on $M$ induced from the Kaehler metric and complex structure $J$ on $M_{n}(c)$.
  The structure vector field $\xi$ is called principal if $A\xi=\alpha\xi$, where A is the
 shape operator of $M$ and $\alpha=\eta(A\xi)$ is a smooth function. A real hypersurface is said to
 be a \textit{Hopf hypersurface} if $\xi$ is principal.

 The classification problem of real hypersurfaces in complex space
 forms is of great importance in Differential Geometry. The study
 of this was initiated by Takagi \cite{T1}, \cite{T2}, who classified all homogenous
 real hypersurfaces in $\mathbb{C}P^{n}$ into six types, which
 are said to be of type $A_{1}$, $A_{2}$, $B$, $C$, $D$ and $E$.
In \cite{CR} Hopf hypersurfaces were considered as tubes over certain
submanifolds in $\mathbb{C}P^{n}$. In \cite{K} the local classification
theorem for Hopf hypersurfaces with constant principal curvatures
in $\mathbb{C}P^{n}$ was given. In the case of complex hyperbolic
space $\mathbb{C}H^{n}$, the classification  theorem for Hopf
hypersurfaces with constant principal curvatures was given by
Berndt \cite{Ber}.

Okumura \cite{Ok}, in $\mathbb{C}P^{n}$, and Montiel and Romero \cite{MR}, in
$\mathbb{C}H^{n}$, gave the classification of real hypersurfaces
satisfying relation $A\varphi=\varphi A$.

\begin{theorem}
Let M be a real hypersurface of $M_{n}(c)$ , $n\geq2$ ($c\neq0$). If it satisfies
$A\varphi-\varphi A=0$, then M is locally congruent to one of the
following hypersurfaces:
 \begin{itemize}
    \item  In case $\mathbb{C}P^{n}$\\
    $(A_{1})$   a geodesic hypersphere of radius r , where
    $0<r<\frac{\pi}{2}$,\\
    $(A_{2})$  a tube of radius r over a totally geodesic
    $\mathbb{C}P^{k}$,$(1\leq k\leq n-2)$, where $0<r<\frac{\pi}{2}.$
    \item In case $\mathbb{C}H^{n}$\\
    $(A_{0})$   a horosphere in $ \mathbb{C}H^{n}$, i.e a Montiel tube,\\
    $(A_{1})$  a geodesic hypersphere or a tube over a hyperplane $\mathbb{C}H^{n-1}$,\\
    $(A_{2}) $  a tube over a totally geodesic $\mathbb{C}H^{k}$ $(1\leq k\leq n-2)$.
  \end{itemize}
\end{theorem}

The Jacobi operator with respect to \emph{X} on M is defined by
 $R(\cdot, X)X$, where R is the Riemmanian curvature of M. For
 $X=\xi$ the Jacobi operator is called structure Jacobi operator
 and is denoted by $l=R(\cdot, \xi)\xi$. It has a
 fundamental role in almost contact manifolds.
 Many differential geometers have studied real hypersurfaces in terms of the structure Jacobi operator.

 The study of real hypersurfaces whose structure Jacobi operator satisfies conditions concerned to the parallelness of it
is a problem of great importance. In \cite{OPS} the nonexistence
of real hypersurfaces in nonflat complex space form with parallel
structure Jacobi operator ($\nabla l=0$) was proved. In \cite{PSaSuh} a
weaker condition ($\mathbb{D}$-parallelness), that is
$\nabla_{X}l=0$ for any vector field $X$ orthogonal to $\xi$, was
studied and it was proved the nonexistence of such real hypersurfaces
in case of $\mathbb{C}P^{n}$ ($n\geq3$). The $\xi$-parallelness of
structure Jacobi operator in combination with other conditions was
another problem that was studied by many other authors such as Ki, Perez, Santos, Suh (\cite{KPSaSuh}).

A tensor field $P$ of type (1, \emph{s}) is said to be
\emph{semi-parallel} if $R\cdot P=0$, where $R$ acts on $P$ as a derivation.

More generally, it is said to be
\emph{pseudo-parallel} if there exists a function \emph{L} such
that $$R\cdot P=L \{(X\wedge Y)\cdot P\},$$ where
$(X\wedge Y)Z=g ( Y, Z )X-g( Z, X)Y$. If $\emph{L}\neq0$, then
the pseudo-parallel tensor is called \emph{proper}.

A Riemannian manifold $M$ is said to be \emph{semi-symmetric }if
$R\cdot R=0$, where the Riemannian curvature tensor $R$ acts on $R$  as a derivation. Deszcz in \cite{De} introduced the notion of \emph{pseudo-symmetry}. A Riemannian manifold is said to
be \emph{pseudo-symmetric} if there exists a function \emph{L }
such that $R( X, Y )\cdot R=L \{(X\wedge Y)\cdot R\}$. If \emph{L} is a constant then the pseudo-symmetric space is called a \emph{pseudo-symmetric
space of constant type}. Both of these notions were studied in the
case of real hypersurfaces in complex space forms. More precisely,
in \cite{NR2} Niebergall and Ryan proved the non-existence of semi-symmetric Hopf real hypersurfaces
and recently in \cite{CHI} Cho, Hamada and Inoguchi gave the classification of pseudo-symmetric Hopf real
hypersurfaces in $\mathbb{C}P^{2}$ and $\mathbb{C}H ^{2}$.

Recently, in \cite{PS} Perez and Santos proved that there exist no real hypersurfaces in complex projective space
$\mathbb{C}P^{n}$, $n\geq3$, with semi-parallel structure Jacobi operator,
(i.e. $R\cdot l=0$). Cho and Kimura in \cite{CK} generalized the previous work and proved the non-existence of
real hypersurfaces in complex space forms, whose structure Jacobi
operator is semi-parallel.

From the above raises naturally the question:
\begin{center}
"Do there exist real hypersurfaces with pseudo-parallel structure Jacobi operator?"
\end{center}

In this paper, we study real hypersurfaces in $\mathbb{C}P^{2}$
and $\mathbb{C}H ^{2}$ equipped with \emph{pseudo-parallel
structure Jacobi operator}, i.e. the structure Jacobi operator
satisfies the following condition:
$$R( X, Y )\cdot l=L \{(X\wedge Y)\cdot l\},$$
more precisely:
\begin{eqnarray}
R(X, Y)lZ-l( R(X,Y)Z)=L\{(X\wedge Y)lZ-l((X\wedge Y)Z)\},
\end{eqnarray}
with $L\neq0$.

Even though Cho and Kurihara proved in [4] the non-existence of real hypersurfaces in complex space form, whose structure Jacobi operator is semi-parallel, in the present paper we prove the existence of real hypersurfaces, whose structure Jacobi operator is pseudo-parallel and we classify them. More precisely:\\

\begin{pro}
Every real hypersurface M in $\mathbb{C}P^{2}$ or
$\mathbb{C}H^{2}$, equipped with pseudo-parallel structure Jacobi operator is a Hopf hypersurface.\\
In case of $\mathbb{C}P^{2}$, M is locally congruent to:
\begin{itemize}
 \item a geodesic hypersphere of radius r, where $0<r<\frac{\pi}{2}$,
 \item or to a non-homogeneous real hypersurface, which is considered as a tube of radius $\frac{\pi}{4}$ over a holomorphic curve in $\mathbb{C}P^{2}$.
\end{itemize}
In case of $\mathbb{C}H^{2}$, M is locally congruent to:
\begin{itemize}
  \item a horosphere,
  \item or to a geodesic hypersphere,
  \item or to a tube over $\mathbb{C}H^{1}$,
  \item or to a Hopf hypersurface with $\eta(A\xi)=0$ in $\mathbb{C}H^{2}$.
\end{itemize}
\end{pro}

\section{Preliminaries}

\hspace{15pt}Throughout this paper all manifolds, vector fields e.t.c. are assumed to be of class $C^{\infty}$ and all manifolds are assumed to be connected. Furthermore, the real hypersurfaces are supposed to be oriented and without boundary. Let $M$ be a real hypersurface immersed in a nonflat
complex space form $(M_{n}(c),G)$ with almost complex structure $J$ of constant holomorphic sectional curvature $c$. Let $N$ be a unit
normal vector field on $M$ and $\xi=-JN$. For a vector field $X$ tangent to $M$ we can write $JX=\varphi (X)+\eta(X)N$, where
$\varphi X$ and $\eta(X)N$ are the tangential and the normal
component of $JX$ respectively. The Riemannian connection
$\overline{\nabla}$ in $M_{n}(c)$ and $\nabla$ in $M$ are related
for any vector fields $X$, $Y$ on $M$:
$$\overline{\nabla}_{Y}X=\nabla_{Y}X+g(AY,X)N,$$
$$\overline{\nabla}_{X}N=-AX,$$
where g is the Riemannian metric on $M$ induced from G of $M_{n}(c)$ and A is the shape operator of $M$ in $M_{n}(c)$. $M$ has an almost contact metric structure $(\varphi,\xi,\eta)$ induced from $J$ on $M_{n}(c)$ where $\varphi$ is a (1,1) tensor field and $\eta$ a
1-form on $M$ such that ([2])
$$g(\varphi X,Y)=G(JX,Y),\hspace{20pt}\eta(X)=g(X,\xi)=G(JX,N).$$
Then we have
\begin{eqnarray}
\varphi^{2}X=-X+\eta(X)\xi,\hspace{20pt}
\eta\circ\varphi=0,\hspace{20pt} \varphi\xi=0,\hspace{20pt}
\eta(\xi)=1,
\end{eqnarray}
\begin{eqnarray}\hspace{20pt}
g(\varphi X,\varphi
Y)=g(X,Y)-\eta(X)\eta(Y),\hspace{10pt}g(X,\varphi Y)=-g(\varphi
X,Y),
\end{eqnarray}
\begin{eqnarray}
\nabla_{X}\xi=\varphi
AX,\hspace{20pt}(\nabla_{X}\varphi)Y=\eta(Y)AX-g(AX,Y)\xi.
\end{eqnarray}
    Since the ambient space is of constant holomorphic sectional
curvature $c$, the equations of Gauss and Codazzi for any vector
fields $X$, $Y$ , $Z$ on $M$ are respectively given by
\begin{eqnarray}
R(X,Y)Z=\frac{c}{4}[g(Y,Z)X-g(X,Z)Y+g(\varphi Y ,Z)\varphi
X\end{eqnarray} $$-g(\varphi X,Z)\varphi Y-2g(\varphi X,Y)\varphi
Z]+g(AY,Z)AX-g(AX,Z)AY,$$
\begin{eqnarray}
\hspace{10pt}
(\nabla_{X}A)Y-(\nabla_{Y}A)X=\frac{c}{4}[\eta(X)\varphi
Y-\eta(Y)\varphi X-2g(\varphi X,Y)\xi],
\end{eqnarray}
where $R$ denotes the Riemannian curvature tensor on $M$.\\
    Relation (2.4) implies that the structure Jacobi operator $l$ is given by:
\begin{eqnarray}
lX=\frac{c}{4}[X-\eta(X)\xi]+\alpha AX-\eta(AX)A\xi.
\end{eqnarray}

    For every point $P\;\;\epsilon\;\; M$, the tangent space
$T_{P}M$ can be decomposed as following:
$$T_{P}M=span\{\xi\}\oplus \mathbb{D}$$
where $\mathbb{D}=\{X\;\;\epsilon\;\; T_{P}M:\eta(X)=0\}$.
 Due to the above decomposition,the vector field $A\xi$ can be written:
 \begin{eqnarray}
 A\xi=\alpha\xi+\beta U,
 \end{eqnarray}
 where $\beta=|\varphi\nabla_{\xi}\xi|$ and
 $U=-\frac{1}{\beta}\varphi\nabla_{\xi}\xi\;\epsilon\;ker(\eta)$, provided
 that $\beta\neq0$.

\section{Some previous results}

 \hspace{15pt}In the rest of this paper, we use the notion $M_{2}(c)$, $c\neq0$, to denote $\mathbb{C}P^{2}$ or
$\mathbb{C}H^{2}$.

 Let $M$ be a non-Hopf hypersurface in $M_{2}(c)$. Then the following relations holds on every three-dimensional real hypersurface in $M_{2}(c)$.
\begin{lemma}
Let M be a real hypersurface in $M_{2}(c)$. Then the following
relations hold on M:
\begin{eqnarray}
\hspace{-70pt}AU=\gamma U+\delta\varphi U+\beta\xi,\hspace{20pt}
A\varphi U=\delta U+\mu\varphi U,
\end{eqnarray}
\begin{eqnarray}
\nabla_{U}\xi=-\delta U+\gamma\varphi U,\hspace{20pt}
\nabla_{\varphi U}\xi=-\mu U+\delta\varphi U,\hspace{20pt}
\nabla_{\xi}\xi=\beta\varphi U,
\end{eqnarray}
\begin{eqnarray}
\nabla_{U}U=\kappa_{1}\varphi U+\delta\xi,\hspace{20pt}
\nabla_{\varphi U}U=\kappa_{2}\varphi U+\mu\xi,\hspace{20pt}
\nabla_{\xi}U=\kappa_{3}\varphi U,
\end{eqnarray}
\begin{equation}
\nabla_{U}\varphi U=-\kappa_{1}U-\gamma\xi,\hspace{5pt}
\nabla_{\varphi U}\varphi U=-\kappa_{2}U-\delta\xi,\hspace{5pt}
\nabla_{\xi}\varphi U=-\kappa_{3}U-\beta\xi,
\end{equation}
where $\gamma,\delta,\mu,\kappa_{1},\kappa_{2},\kappa_{3}$ are
smooth functions on M.
\end{lemma}
\textbf{Proof:} Let $\{U,\varphi U,\xi\}$ be an orthonormal basis of $M$. Then we
have:
$$\hspace{10pt}AU=\gamma U+\delta\varphi U+\beta\xi\hspace{30pt}A\varphi U=\delta U+\mu\varphi U,$$ where $\gamma,\delta,\mu$ are smooth functions, since $g(AU,\xi)=g(U,A\xi)=\beta$ and
$g(A\varphi U,\xi)=g(\varphi U,A\xi)=0$.

The first relation of (2.3), because of (2.6) and (3.1), for $X=U$,
$X=\varphi U$ and $X=\xi$ implies (3.2).

 From the well known
relation: $Xg(Y,Z)=g(\nabla_{X}Y,Z)+g(Y,\nabla_{X}Z)$ for
$X,Y,Z$ $\epsilon$ $\{\xi,U,\varphi U\}$ we obtain (3.3) and
(3.4), where $\kappa_{1}, \kappa_{2}$ and $\kappa_{3}$ are smooth
functions.\qed
\\

In [7], T.A.Ivey and P.J.Ryan proved the non-existence of real hypersurfaces in $M_{2}(c)$, whose structure Jacobi operator vanishes. In our context, we give a different proof of their Proposition 8 (non-Hopf case) and Lemma 9.
\begin{proposition}
There does not exist real non-flat hypersurface in $M_{2}(c)$, whose structure Jacobi operator vanishes.
\end{proposition}
\textbf{Proof:} Let $M$ be a non-Hopf real hypersurface in $M_{2}(c)$, so the vector field $A\xi$ can be written $A\xi=\alpha\xi+\beta U$ (i.e. $\alpha\beta\neq0$).

Let $\{U,\varphi U,\xi\}$ denote an orthonormal basis of $M$.
Since the structure Jacobi operator of $M$ vanishes, from relation (2.6) for $X=U$ and $X=\varphi U$, we obtain: $AU=(\frac{\beta^{2}}{\alpha}-\frac{c}{4\alpha})U+\beta\xi$ and $A\varphi U=-\frac{c}{4\alpha}\varphi U$. Conversely, if we have a real hypersurface, whose shape operator satisfies the last relations then $l=0$. Relations (3.2), (3.3) and (3.4) because of the latter become respectively:
\begin{eqnarray}
\hspace{-55pt}\nabla_{U}\xi=(\frac{\beta^{2}}{\alpha}-\frac{c}{4\alpha})\varphi U,\hspace{20pt}
\nabla_{\varphi U}\xi=\frac{c}{4\alpha}U,\hspace{20pt}
\nabla_{\xi}\xi=\beta\varphi U,
\end{eqnarray}
\begin{eqnarray}
\hspace{-40pt}\nabla_{U}U=\kappa_{1}\varphi U,\hspace{20pt}
\nabla_{\varphi U}U=\kappa_{2}\varphi U-\frac{c}{4\alpha}\xi,\hspace{20pt}
\nabla_{\xi}U=\kappa_{3}\varphi U,
\end{eqnarray}
\begin{equation}
\hspace{-10pt}\nabla_{U}\varphi U=-\kappa_{1}U-(\frac{\beta^{2}}{\alpha}-\frac{c}{4\alpha})\xi,\hspace{5pt}
\nabla_{\varphi U}\varphi U=-\kappa_{2}U,\hspace{5pt}
\nabla_{\xi}\varphi U=-\kappa_{3}U-\beta\xi,
\end{equation}
where $\kappa_{1},\kappa_{2},\kappa_{3}$ are
smooth functions on M.

   On $M$ the Codazzi equation for $X$, $Y$ $\epsilon$ $\{U,\varphi
U\,\xi\}$, because of (3.5), (3.6) and (3.7) yields:
\begin{eqnarray}
U\beta&=&\beta\kappa_{2}(\frac{4\beta^{2}}{c}+1),\\
\frac{\beta^{2}\kappa_{3}}{\alpha}&=&\beta\kappa_{1}+\frac{c}{4\alpha}(\frac{\beta^{2}}{\alpha}-\frac{c}{4\alpha}),\\
U\alpha&=&\xi\beta\hspace{5pt}=\hspace{5pt} \frac{4\alpha\beta^{2}\kappa_{2}}{c},\\
\xi\alpha&=&\frac{4\alpha^{2}\beta\kappa_{2}}{c},\\
(\varphi U)\alpha&=&\beta(\alpha+\kappa_{3}+\frac{3c}{4\alpha}),\\
(\varphi
U)\beta&=&\beta^{2}+\beta\kappa_{1}+\frac{c}{2\alpha}(\frac{\beta^{2}}{\alpha}-\frac{c}{4\alpha}),\\
(\varphi
U)(\frac{\beta^{2}}{\alpha}-\frac{c}{4\alpha})&=&\beta(\frac{\beta^{2}}{\alpha}+\frac{\beta\kappa_{1}}{\alpha}-\frac{3c}{4\alpha}).
\end{eqnarray}

The Riemannian curvature on $M$ satisfies (2.4) and on the other
hand is given by the relation $R(X, Y)
Z=\nabla_{X}\nabla_{Y}Z-\nabla_{Y}\nabla_{X}Z-\nabla_{[X,Y]}Z$. The combination of these two relations implies:
\begin{eqnarray}
U\kappa_{3}-\xi\kappa_{1}&=&\kappa_{2}(\frac{\beta^{2}}{\alpha}-\frac{c}{4\alpha}-\kappa_{3}),\\
(\varphi
U)\kappa_{3}-\xi\kappa_{2}&=&\kappa_{1}(\kappa_{3}+\frac{c}{4\alpha})+\beta(\kappa_{3}-\frac{c}{2\alpha}).
\end{eqnarray}
Relation (3.14), because of (3.9), (3.12) and (3.13), yields:
\begin{eqnarray}
\kappa_{3}=-4\alpha,
\end{eqnarray}
  and so relation (3.9) becomes:
\begin{eqnarray}
\beta\kappa_{1}=\frac{c}{4\alpha}(\frac{c}{4\alpha}-\frac{\beta^{2}}{\alpha})-4\beta^{2}.
\end{eqnarray}

Differentiating the relations (3.17) and (3.18) with respect to U
and $\xi$ respectively and substituting in (3.15) and due to
(3.10), (3.11) and (3.17)  we obtain:
\begin{eqnarray}
\kappa_{2}(c-2\beta^{2}-4\alpha^{2})=0.
\end{eqnarray}

Owing to (3.19), we consider $M_{1}$ the open subset of points
$P\;\epsilon\; M$, where $\kappa_{2}\neq0$ in a neighborhood of
every $P$. Due to (3.19) we obtain: $2\beta^{2}+4\alpha^{2}=c$ on $M_{1}$.
Differentiation of the last relation along $\xi$ and taking into
account (3.10), (3.11) and $2\beta^{2}+4\alpha^{2}=c$  yields:
$c=0$, which is a contradiction. Therefore, $M_{1}$ is empty. Thus, $\kappa_{2}=0$ on $M$ and relations (3.8), (3.10) and (3.11) become:
$$U\alpha=U\beta=\xi\alpha=\xi\beta=0.$$
Using the above relations we obtain:
$$\hspace{-10pt}
[U, \xi]\alpha=U\xi\alpha-\xi U\alpha=0,$$
$$[U,\xi]\alpha=(\nabla_{U}\xi-\nabla_{\xi}U)\alpha=\frac{1}{4\alpha}(4\beta^{2}+16\alpha^{2}-c)(\varphi
U)\alpha.$$ Combining the last two relations we have:
\begin{eqnarray}
(4\beta^{2}+16\alpha^{2}-c)(\varphi U)\alpha=0.
\end{eqnarray}

Let $M_{2}$ be the set of points $P\;\epsilon\; M$, for which
there exists a neighborhood of every P such that $(\varphi
U)\alpha\neq0$. So in $M_{2}$ from (3.20) we have:
$16\alpha^{2}+4\beta^{2}=c$. Differentiating the last relation
with respect to $\varphi U$ and taking into account (3.12),
(3.13), (3.17), (3.18) and $16\alpha^{2}+4\beta^{2}=c$, we obtain: $4\alpha^{2}+\beta^{2}=0$,
which is impossible. So $M_{2}$ is empty. Hence, on $M$ we have
$(\varphi U)\alpha=0$. Then, relations (3.12), (3.17) and (3.18) imply: $c=4\alpha^{2}$ and
$\beta\kappa_{1}=\alpha^{2}-5\beta^{2}$. On the other hand from relation (3.16), because of (3.17) we obtain:
$\kappa_{1}=-2\beta$. Substitution of $\kappa_{1}$ in
$\beta\kappa_{1}=\alpha^{2}-5\beta^{2}$ yields:
$3\beta^{2}=\alpha^{2}$. Taking the covariant derivative along $\varphi U$ of $3\beta^{2}=\alpha^{2}$, because of (3.13), we conclude: $\beta=0$,
which is a contradiction.

Suppose that $A\xi=\beta\xi$ (i.e. $\alpha=0$ and $\beta\neq0$). Since the structure Jacobi operator of $M$ vanishes, from relation (2.6) for $X=\varphi U$, we obtain: $c=0$, which is impossilbe.

Hence, there do not exist non-Hopf hypersurfaces with $l=0$. Using this and the Hopf case (\cite{IR}), we complete the proof of the present Proposition.
\qed

\section{Auxiliary Relations}
\hspace{15pt} If $M$ is a real hypersurface in $M_{2}(c)$, we  consider the open subset $\mathcal{N}$ of $M$ such that:
$$\mathcal{N}=\{P\;\;\epsilon\;\;M:\;\beta\neq0,\;\;\mbox{in neighborhood of P}\}.$$
Furthermore, we consider $\mathcal{V}$, $\Omega$ open subsets of $\mathcal{N}$ such that:
$$\mathcal{V}=\{P\;\;\epsilon\;\;\mathcal{N}:\alpha=0,\;\;in\;\;a\;\;neighborhood\;\;of\;\;P\},$$
$$\Omega=\{P\;\;\epsilon\;\;\mathcal{N}:\alpha\neq0,\;\;in\;\;a\;\;neighborhood\;\;of\;\;P\},$$
where $\mathcal{V}\cup\Omega$ is open and dense in the closure of $\mathcal{N}$.

\begin{lemma}
Let M be a real hypersurface in $M_{2}(c)$, equipped with
pseudo-parallel structure Jacobi operator. Then $\mathcal{V}$ is empty.
\end{lemma}
\textbf{Proof:} Let $\{U,\varphi U,\xi\}$ be a local orthonormal basis on $\mathcal{V}$. The
relation (2.7) takes the form $A\xi=\beta U$ and we consider:
\begin{eqnarray}
AU=\gamma' U+\delta'\varphi U+\beta\xi,\hspace{20pt} A\varphi U=\delta'U+\mu'\varphi U,
\end{eqnarray}
since $g(AU,\xi)=g(U,A\xi)=\beta$, $g(A\varphi U,\xi)=g(\varphi
U,A\xi)=0$ and $\gamma',\delta',\mu'$ are smooth
functions.\\
From (2.6) for $X=U$ and $X=\varphi U$, taking into account (4.1), we obtain:
\begin{eqnarray}
l\varphi U=\frac{c}{4}\varphi U\hspace{20pt}lU=(\frac{c}{4}-\beta^{2})U.
\end{eqnarray}

Relation (1.1) for $X=U$, $Y=\xi$ and $Z=\varphi U,$ because of (2.4), (4.1) and (4.2) yields: $\delta'=0$, since $\beta\neq0$.

Furthermore, relation (1.1) for $X=U$ and $Y=Z=\varphi U$, owing to (2.4), (4.1), (4.2) and $\delta'=0$ implies:
\begin{eqnarray}
\mu'=0 \hspace{20pt}c=L,
\end{eqnarray}
and  for $X=\xi$ and $Y=Z=\varphi U$, because of (4.3), gives: $c=0$, which is a
contradiction. Therefore, $\mathcal{V}$ is empty.
\qed
\\

In what follows we work on $\Omega$, where $\alpha\neq0$ and $\beta\neq0$.

 By using (2.6) and relations (3.1) we obtain:
\begin{eqnarray}
\hspace{10pt}
lU=(\frac{c}{4}+\alpha\gamma-\beta^{2})U+\alpha\delta\varphi
U\hspace{20pt}l\varphi U=\alpha\delta
U+(\alpha\mu+\frac{c}{4})\varphi U
\end{eqnarray}
The relation (1.1) because of (2.4), (3.1) and (4.4), implies:
\begin{eqnarray}
\delta&=&0,\;\;\mbox{for $X=U$,\;\;$Y=\xi$\;\;and\;\;$Z$=$\varphi U$},
\end{eqnarray}
and additional due to (4.5) yields:
\begin{eqnarray}
 \mu(\alpha\mu+\frac{c}{4})&=&0,\;\;\mbox{for $X=U$,\;\;$Y=\varphi U$\;\;and\;\;$Z=\xi$}.
\end{eqnarray}

Owing to (4.6), we consider $\Omega_{1}$ the open subset of $\Omega$, such that:
$$\Omega_{1}=\{P\;\;\epsilon\;\;\Omega:\mu\neq-\frac{c}{4\alpha},\;\;in\;\;a\;\;neighborhood\;\;of\;\;P\}.$$
Therefore, in $\Omega_{1}$ from (4.6) we have: $\mu=0$ .

\begin{lemma}
Let M be a real hypersurface in $M_{2}(c)$, equipped with
pseudo-parallel structure Jacobi operator. Then $\Omega_{1}$ is empty.
\end{lemma}
\textbf{Proof:} In $\Omega_{1}$, relation (1.1) for $X=U$, $Y=\varphi U$ and $Z=U$, because of (2.4), (3.1), (4.4) and (4.5) yields:
\begin{eqnarray}
(\beta^{2}-\alpha\gamma)(c-L)=0.
\end{eqnarray}
Due to (4.7), we consider the open subset $\Omega_{11}$ of $\Omega_{1}$, such that:
$$\Omega_{11}=\{P\;\;\epsilon\;\;\Omega_{1}:c\neq L,\;\;in\;\;a\;\;neighborhood\;\;of\;\;P\}.$$
So in $\Omega_{11}$, we obtain: $\gamma=\frac{\beta^{2}}{\alpha}$.

In $\Omega_{11}$, the relation (2.5), because of Lemma 3.1 and (4.5), yields:
\begin{eqnarray}
\frac{\beta^{2}\kappa_{3}}{\alpha}&=&\beta\kappa_{1}+\frac{c}{4},\;\;\mbox{for $X=U$\;\;and\;\;$Y=\xi$}\\
(\varphi U)\alpha&=&\beta(\alpha+\kappa_{3}),\;\;\mbox{for $X=\varphi$U\;\;and\;\;$Y=\xi$}\\
(\varphi
U)\beta&=&\beta^{2}+\beta\kappa_{1}+\frac{c}{2},\;\;\mbox{for $X=\varphi$U\;\;and\;\;$Y=\xi$}\\
(\varphi
U)\frac{\beta^{2}}{\alpha}&=&\frac{\beta^{2}}{\alpha}(\kappa_{1}+\beta),\;\;\mbox{for $X=U$\;\;and\;\;$Y=\varphi U$}.
\end{eqnarray}
Substituting in (4.11) the relations (4.9), (4.10) and taking
into account (4.8) we obtain: $\frac{3c\beta}{4\alpha}=0$, which is
a contradiction. Therefore, $\Omega_{11}$ is empty and  $L=c$ in
$\Omega_{1}$.

In $\Omega_{1}$, relation (1.1) for $X=\xi$ and $Y=Z=\varphi U$, because of (2.4), (3.1) and (4.4) implies:
$c=0$, which is impossible. Therefore, $\Omega_{1}$ is empty. \qed
\\

From Lemma 4.1, we conclude that $\mu=-\frac{c}{4\alpha}$ in
$\Omega$.
\begin{lemma}
Let M be a real hypersurface in $M_{2}(c)$, equipped with pseudo-parallel structure Jacobi operator. Then $\Omega$ is empty.
\end{lemma}
\textbf{Proof:} In $\Omega$, relation (1.1) for $X=\varphi U$, $Y=\xi$ and $Z=U$, due to (2.4), (3.1), (4.4) and (4.5) yields:
$\gamma=\frac{\beta^{2}}{\alpha}-\frac{c}{4\alpha}$. Owing to $\mu=-\frac{c}{4\alpha}$ and $\gamma=\frac{\beta^{2}}{\alpha}-\frac{c}{4\alpha}$ and (4.5), relation (4.4) implies: $lU=l\varphi U=0$ and since $l\xi=0$, we obtain that the structure Jacobi operator vanishes in $\Omega$. Due to Proposition 3.2, we conclude that $\Omega$ is empty.
\qed
\\

From Lemmas 4.1 and 4.3, we conclude that $\mathcal{N}$ is empty and we lead to the following result:
\begin{proposition}
Every real hypersurface in $M_{2}(c)$, equipped with pseudo-parallel structure Jacobi operator, is a Hopf hypersurface.
\end{proposition}

\section{Proof of Main Theorem}
Since $M$ is a Hopf hypersurface, due to Theorem 2.1 (\cite{NR1})  we have that $\alpha$ is a constant. We consider a unit vector field $e$ $\epsilon$ $\mathbb{D}$, such that $Ae=\lambda e$, then $A\varphi e=\nu\varphi e$ at some point $P$ $\epsilon$ $M$, where $\{ e, \varphi e, \xi\}$ is a local orthonormal basis. Then the following relation holds on $M$, (Corollary 2.3 \cite{NR1}):
\begin{eqnarray}
\lambda\nu=\frac{\alpha}{2}(\lambda+\nu)+\frac{c}{4}.
\end{eqnarray}
The relation (2.6) implies:
\begin{eqnarray}
le=(\frac{c}{4}+\alpha\lambda)e\;\;\;\mbox{and}\;\;\;l\varphi e=(\frac{c}{4}+\alpha\nu)\varphi e.
\end{eqnarray}
Relation (1.1) for $X=e$ and $Y=Z=\varphi e$, because of (2.4) and (5.2) yields:
\begin{eqnarray}
\alpha(c+\lambda\nu-L)(\nu-\lambda)=0.
\end{eqnarray}
Relation (1.1) for $X=Z=e$, $Y=\xi$ and for $X=Z=\varphi e$, $Y=\xi$, because of (2.4) and (5.2) implies respectively:
\begin{eqnarray}
(\frac{c}{4}+\alpha\lambda)(L-\alpha\lambda-\frac{c}{4})=0,\\
(\frac{c}{4}+\alpha\nu)(L-\alpha\nu-\frac{c}{4})=0.
\end{eqnarray}
Because of (5.3), we consider $\mathcal{M}_{1}$ the open subset of $M$, such that:
$$\mathcal{M}_{1}=\{P\;\;\epsilon\;\;M\;:\alpha(\nu-\lambda)\neq0\;\;in\;\;a\;\;neighborhood\;\;of\;\;P\}.$$
So in $\mathcal{M}_{1}$, we have: $L=c+\lambda\nu$.
\begin{proposition}
Let $M$ be a real Hopf hypersurface in $M_{2}(c)$, equipped with pseudo-parallel structure Jacobi operator. Then $\mathcal{M}_{1}$ is empty.
\end{proposition}
\textbf{Proof:} Because of (5.4), we consider $\mathcal{M}_{11}$  the open subset of $\mathcal{M}_{1}$, such that:
 $$\mathcal{M}_{11}=\{P\;\epsilon\;\mathcal{M}_{1}:\;L\neq\alpha\lambda+\frac{c}{4},\;in\;a\;neighborhood\;of\;P\}.$$
In $\mathcal{M}_{11}$ relations  (5.4) and (5.5) imply: $\lambda=-\frac{c}{4\alpha}$ and $L=\alpha\nu+\frac{c}{4}$, respectively since $\lambda\neq\nu$. Using the last two relations and because of $L=c+\lambda\nu$ and (5.1), we obtain:
\begin{eqnarray}
\lambda=\frac{4\alpha}{7},\hspace{20pt}\nu=-4\alpha,\hspace{20pt}c=-\frac{16\alpha^{2}}{7}.
\end{eqnarray}
Because of (5.6), we have $c<0$ and three distinct constant eigenvalues. So the only case is real hypersurface of type B in $\mathbb{C}H^{2}$.
Substitution of the eigenvalues of type B real hypersurfaces (see \cite{Ber}) in (5.6), leads to a contradiction. So $\mathcal{M}_{11}=\emptyset$. Consequently, in $\mathcal{M}_{1}$ the relation $L=\alpha\lambda+\frac{c}{4}$ holds and because of (5.5), we lead to: $\nu=-\frac{c}{4\alpha}$, since $\lambda\neq\nu$. Following the same method as above, we obtain a contradiction and this completes the proof of the Proposition.\qed
\\

Thus from Proposition 5.1, we conclude that $\alpha(\nu-\lambda)=0$ at any point $P$ $\epsilon$ $M$. Thus locally either $\alpha=0$ or $\nu=\lambda$.

If $\alpha=0$  in case of $\mathbb{C}P^{2}$, $M$ is locally congruent to a tube of radius $r=\frac{\pi}{4}$ over a holomorphic curve in $\mathbb{C}P^{2}$, if $\lambda\neq\nu$ or to a geodesic hypersphere of radius $r=\frac{\pi}{4}$, if $\lambda=\nu$, (see \cite{CR}), and in case of $\mathbb{C}H^{2}$, $M$ is a Hopf hypersurface with $A\xi=0$.

If $\alpha\neq0$, we have: $\lambda=\nu$. Then $Ae=\lambda e$ and $A\varphi e=\lambda\varphi e$, therefore we obtain:  $$(A\varphi-\varphi A)X=0,\;\;\forall\;\;X\;\;\epsilon\;\;TM.$$
From the above relation Theorem 1.1 holds and this completes the proof of Main Theorem.

\section*{Acknowledgements}
 The authors thank Prof. F. Gouli-Andreou for her comments on improving the proof of main theorem.\\
 The first author is granted by the Foundation Alexandros S. Onasis. Grant Nr: G ZF 044/2009-2010

\end{document}